\input amstex
\input epsf
\documentstyle{amsppt}
\magnification 1200
\vcorrection{-12mm}

\chardef\und`\_       

\def\refBGGMTW     {1}      
\def\refPlantriPap {2}      
\def\refPlantri    {3}      
\def\refC          {4}      
\def\refD          {5}      
\def\refN          {6}      
\def\refPPF        {7}

\def\figBNS       {1}
\def\figBNRP      {2}
\def\figN         {3}
\def\figMWheel    {4}
\def\figGener     {5}
\def\figAnnuli    {6}
\def\figPropS     {7}

\def\figContract  {11}
\def\figRemoval   {12}
\def\figB         {13}
\def\figGfive     {14}
\def\figA         {15}
\def\figPPF       {16}

\def\thNakamoto   {1}
\def\thMain       {2}
\def\remDegFour   {3}
\def\propOne      {4}
\def\propLoopFree {5}
\def\lemZero      {6}
\def\propS        {7}
\def\thRP         {8}
\def\thUniq       {9}
\def\lemMin       {10}
\def\lemRemove    {11}
\def\thMin        {12}

\def\sectDef      {2}
\def\sectDefGS       {\sectDef.2}

\def\sectCover    {4}
\def\sectGener    {5}
\def\sectComput   {6}

\def\eqDeg     {1}
\def\eqCommand {2}

\def\RP{\Bbb{RP}}
\def\R{\Bbb R}
\def\Z{\Bbb Z}
\def\id{\operatorname{id}}
\def\trace{\operatorname{trace}}
\def\bip{\operatorname{bip}}

\topmatter
\title     Basic nets in the projective plane
\endtitle
\author    S.~Yu.~Orevkov
\endauthor
\abstract
The notion of basic net (called also basic polyhedron) on $S^2$ plays a central role
in Conway's approach to enumeration of knots and links in $S^3$.
Drobotukhina applied this approach for links in $\RP^3$ using basic nets on $\RP^2$.
By a result of Nakamoto, all basic nets on $S^2$ can be obtained
from a very explicit family of minimal basic nets (the nets $(2\times n)^*$, $n\ge3$,
in Conway's notation) by two local transformations.
We prove a similar result for basic nets in $\RP^2$.

We prove also that a graph on $\RP^2$ is uniquely determined by its
pull-back on $S^2$ (the proof is based on Lefschetz fix point theorem).
\endabstract
\address
IMT, Universit\'e Paul Sabatier (Toulouse-3)
\endaddress
\email orevkov\@math.ups-tlse.fr
\endemail
\address
Steklov Math. Institute, Moscow
\endaddress
\endtopmatter

\document

\head 1. Introduction and statement of main results
\endhead

In this paper, a {\it surface} is a smooth compact 2-manifold without boundary. 
A {\it net} on a surface $F$ is the image of a generic immersion of several
circles.

A net $\Gamma$ is called {\it irreducible} if for any embedded circle
transversally intersecting $\Gamma$ at most at two points
and dividing $F$ into two components $F_1$ and $F_2$, the following
condition holds: one of $F_1$, $F_2$
is a disk whose intersection with $\Gamma$ is either a simple arc or empty.

A net $\Gamma$ on $F$ is called {\it basic} if it is 
irreducible and none of the components of $F\setminus\Gamma$
is a digon whose corners are at two distinct vertices.

Basic nets on $S^2$ (called in [\refC] basic polyhedra) were introduced by
Conway [\refC] as a tool for classification of links in $S^3$. Drobotukhina [\refD]
applied Conway's approach for links in $\RP^3$ using basic nets on $\RP^2$

Basic nets on $S^2$ (resp. on $\RP^2$) with $\le 11$ (resp. with $\le 8$) crossings
are shown in Figure \figBNS\ (resp. in Figure \figBNRP). In all pictures, we represent
$\RP^2$ as a disk whose opposite boundary points are supposed to be identified.
%
The nets $g_i^n$ in Figure \figBNRP\ are denoted as in [\refD] for $n\le6$ (except
that $g^6_3$ is missing in [\refD]; note that the corresponding alternative link
in $\RP^3$ also is missing in [\refD]). For $n\ge7$ we number them in
the order they are produced by {\tt plantri} program [\refPlantri] (see \S\sectComput).

An algorithm to generate all basic nets on $S^2$ with a given number of crossings is
obtained in [\refN], improved in [\refBGGMTW], and implemented in [\refPlantri].
The main purpose of the present paper is to extend these results to $\RP^2$.

We prove also Theorem \thUniq\ which could be of independent interest. It states that
any cellular graph on $\RP^2$ (i.~e. a graph whose complement is a union of open disks)
is uniquely determined by its covering on $S^2$. The proof is based on the Lefschetz
fixed point theorem.

\midinsert
\centerline{\epsfxsize=120mm\epsfbox{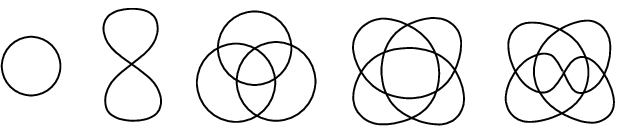}}
\vskip-4mm
\centerline{
    $0^*$\hskip 14mm
    $1^*$\hskip 17mm
    $6^*$\hskip 25mm
    $8^*$\hskip 26mm
    $9^*$\hskip 25mm
}
\vskip2mm
\centerline{\epsfxsize=120mm\epsfbox{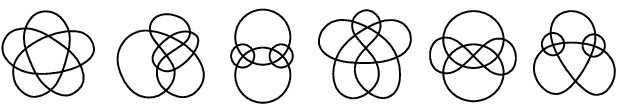}}
\vskip-3mm
\centerline{
    $10^*$\hskip 16mm
    $10^{**}$\hskip 14mm
    $10^{***}$\hskip 12mm
    $11^*$\hskip 16mm
    $11^{**}$\hskip 14mm
    $11^{***}$\hskip 14mm
}
\botcaption{Figure \figBNS}
Basic nets on $S^2$ with $\le 11$ nodes
\endcaption
\endinsert

\midinsert
\centerline{\epsfxsize=125mm\epsfbox{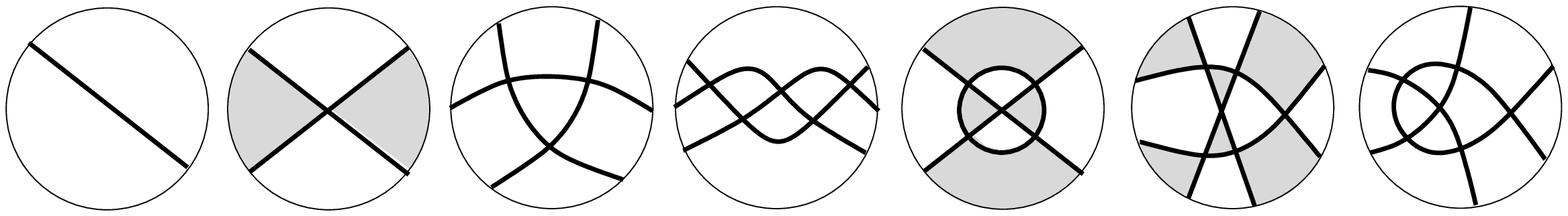}}
\vskip-2mm
\centerline{
    $g^0$  \hskip 14mm
    $g^1$  \hskip 14mm
    $g^3$  \hskip 14mm
    $g^5_1$\hskip 14mm
    $g^5_2$\hskip 14mm
    $g^6_1$\hskip 14mm
    $g^6_2$\hskip 14mm
}
\vskip2mm
\centerline{\epsfxsize=125mm\epsfbox{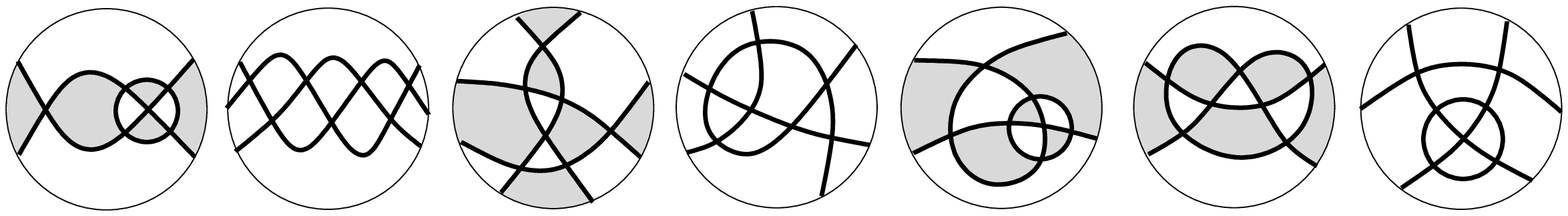}}
\vskip-2mm
\centerline{
    $g^6_3$  \hskip 14mm
    $g^7_1$  \hskip 14mm
    $g^7_2$  \hskip 14mm
    $g^7_3$\hskip 14mm
    $g^7_4$\hskip 14mm
    $g^7_5$\hskip 14mm
    $g^7_6$\hskip 14mm
}
\vskip2mm
\centerline{\epsfxsize=108mm\epsfbox{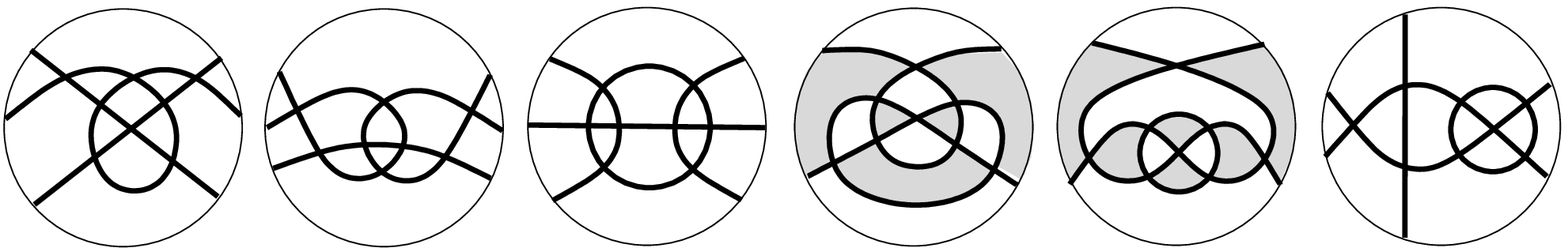}}
\vskip-2mm
\centerline{
    $g^8_1$  \hskip 14mm
    $g^8_2$  \hskip 14mm
    $g^8_3$  \hskip 14mm
    $g^8_4$\hskip 14mm
    $g^8_5$\hskip 14mm
    $g^8_6$\hskip 14mm
}
\vskip2mm
\centerline{\epsfxsize=108mm\epsfbox{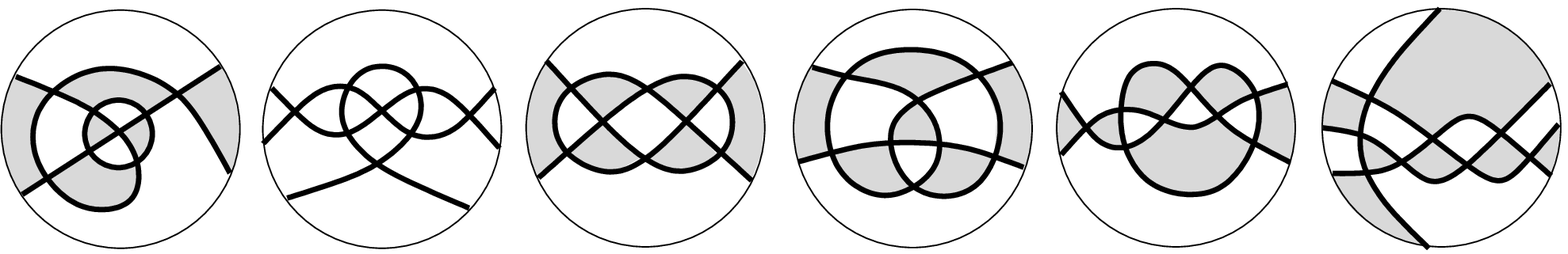}}
\vskip-2mm
\centerline{
    $g^8_7$  \hskip 14mm
    $g^8_8$  \hskip 14mm
    $g^8_9$  \hskip 14mm
    $g^8_{10}$\hskip 13mm
    $g^8_{11}$\hskip 13mm
    $g^8_{12}$\hskip 13mm
}
\botcaption{Figure \figBNRP}
Basic nets on $\RP^2$ with $\le 8$ nodes
\endcaption
\endinsert

%
Given a basic net $\Gamma$ of a surface $F$, one can obtain another net using
the following transformations: 

\roster
\item {\it Face splitting.} Suppose that one of the faces of
$\Gamma$ is an $n$-gon $f$, $n\ge4$. Let $\gamma$ be a simple arc
inside $f$ which connects two non-consecutive sides of $f$ represented by two
distinct edges of $\Gamma$ (see \S\sectDefGS\ for a definition of faces, sides and edges).
Then a neighbourhood of $\gamma$ is replaced as in Figure \figN.1.

\item {\it Vertex surrounding.} A neighbourhood of a vertex of $\Gamma$
is replaced as in Figure \figN.2
\endroster
We say that a face splitting is {\it special} if an $n$-gon splits
into an $(n-1)$-gon and a triangle.



\midinsert
\centerline{
\hskip0pt
\hbox{\vbox{\hsize=50mm
  \epsfxsize=30mm
  \centerline{\epsfbox{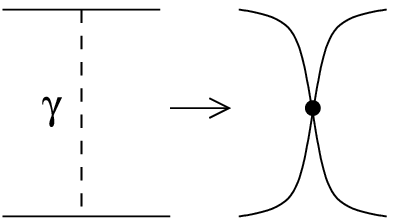}}
  \botcaption{Figure \figN.1} Face splitting
  \endcaption
}}
\hskip 10mm
\hbox{\vbox{\hsize=50mm
  \epsfxsize=40mm
  \centerline{\epsfbox{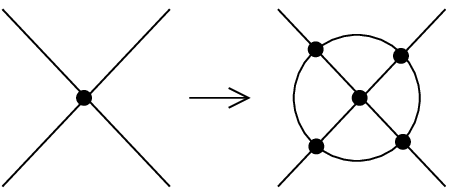}}
  \botcaption{Figure \figN.2} Vertex surrounding
  \endcaption
}}
}
\endinsert

The nets $6^*$, $8^*$, and $10^*$ (Figure \figBNS) are the first three members
of a series denoted in
[\refC] by $(2\times n)^*$, $n\ge 3$. The net $(2\times n)^*$ (up to homeomorphism of $S^2$)
is the union of a regular $n$-gon with its inscribed and circumscribed circles.
In [\refN], [\refBGGMTW], the dual graph of $(2\times n)^*$ is called
{\it pseudo-double wheel} and is denoted by $W_n$.

A result of Nakamoto [\refN; Theorem 1] improved in [\refBGGMTW; Theorem 2]
can be reformulated as follows.

\proclaim{ Theorem 1 } All basic nets on $S^2$ except $0^*$ and $1^*$
can be obtained from $(2\times n)^*$, $n\ge 3$, by successive special face splittings
and vertex surroundings.
\endproclaim


In this paper, we generalize Theorem 1 for basic nets in $\RP^2$.
For an odd $n\ge3$, let $\overline{(2\times n)^*}$ be the net in $\RP^2$ whose
double covering is the net $(2\times n)^*$ on $S^2$. Due to Theorem \thUniq\ in \S\sectCover,
$\overline{(2\times n)^*}$ is uniquely determined by this condition.
It can be described also as follows. Let $P$ be a regular $n$-gon inscribed in
a circle $S$ which bounds a disk $D$. Then, up to homeomorphism, $\overline{(2\times n)^*}$
is the image of $P\cup S$ on the projective plane obtained from $D$ by identifying the
opposite boundary points. For $n=3,5,7$, these are the nets $g^3,g^5_1,g^7_1$ in Figure \figBNRP.
The dual graph of the net $\overline{(2\times n)^*}$ is called in [\refN] {\it M\"obius wheel}
and is denoted by $\widetilde W_n$ (see Figure \figMWheel).

\midinsert
\centerline{\epsfxsize=60mm\epsfbox{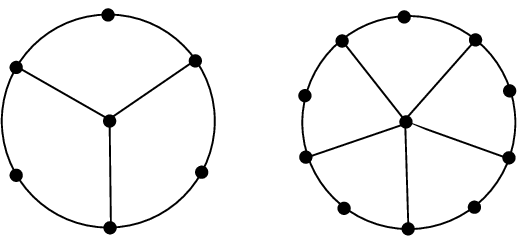}}
\botcaption{Figure \figMWheel}
M\"obius wheels $\widetilde W_3$ and $\widetilde W_5$ in $\RP^2$
\endcaption
\endinsert

We say that a net in $\RP^2$ is {\it homologically trivial} (resp.
{\it homologically non-trivial\/}) if it represents zero (resp. non-zero)
homology class in $H_1(\RP^2;\Z_2)$; in this case the
dual graph is bipartite (resp. non-bipartite). In Figure \figBNRP\
we use the chess-board coloring for the homologically trivial nets.
It is easy to see that the face splittings and vertex surroundings
do not change the homology class.

\proclaim{ Theorem \thMain }
(a). All homologically trivial basic nets on $\RP^2$
can be obtained from $g^1$ {\rm(}see Figure \figBNRP\/{\rm)}
by successive special face splittings and vertex surroundings.

(b). All homologically non-trivial basic nets on $\RP^2$ except $g^0$
can be obtained from $\overline{(2\times n)^*}$ with odd $n\ge 3$
by successive special face splittings and vertex surroundings.
\endproclaim

We prove this theorem in \S\sectGener.

In Figure \figGener, we show all the possible special face splittings and vertex surroundings
on the basic nets in $\RP^2$ with $\le8$ crossings. The number of different
special face splittings which produce the same result is indicated in parentheses near each arrow.
Thus, the list in Figure \figBNRP\ is exhaustive by Theorem \thMain.

\midinsert
\epsfxsize=100mm
\centerline{\epsfbox{ 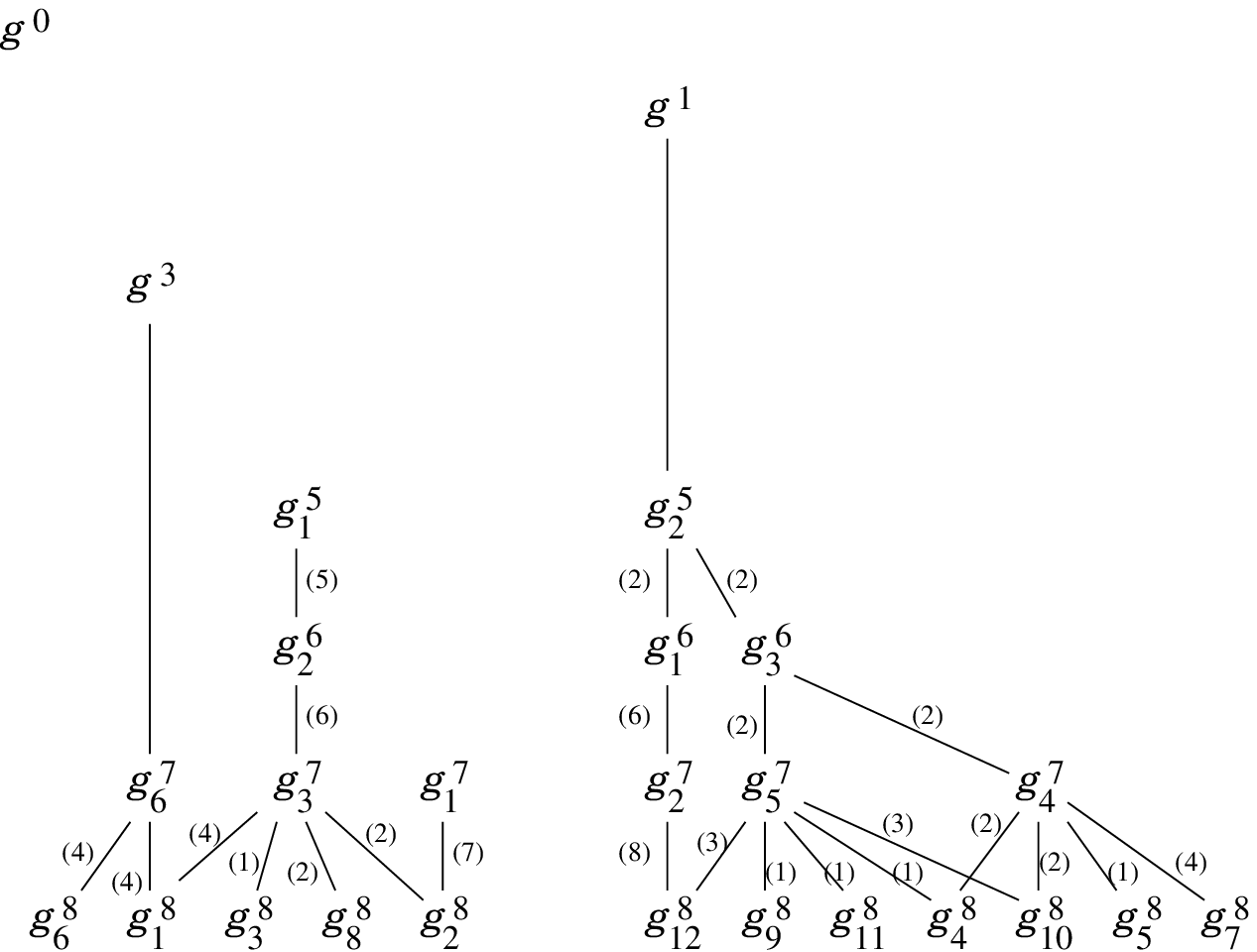 }}
\botcaption{ Figure \figGener } Generating basic nets on $\RP^2$ with $\le8$ crossings
\endcaption
\endinsert

%
Theorem 1 provides an algorithm to generate all basic nets in $S^2$. This
algorithm is efficiently implemented in the program {\tt plantri}
[\refPlantriPap, \refPlantri].
Theorem \thMain\ provides a similar algorithm for $\RP^2$ but it is not implemented
yet. Instead, we used {\tt plantri} with a simple additional filter [\refPPF]
to generate all basic nets on $\RP^2$ up to 18 crossings, see \S\sectComput\
for more details.

\subhead Acknowledgment \endsubhead I am grateful to the referee for indicating
some mistakes in the first version of this paper.


\head 2. Definitions
\endhead
\subhead 2.1. Graphs
\endsubhead
We shall use the following terminology.
A {\it graph} $G$ is a triple $(V,E,\partial)$ where $V=V(G)$ and $E=E(G)$
are two sets whose elements are called {\it vertices} and {\it edges}
respectively, and $\partial$ is a mapping from $E(G)$ to the set of unordered 
pairs of vertices. If $\partial(e)=\{a,b\}$, then $a$ and $b$ are called the
{\it ends} of $e$. A graph $G$ is called {\it finite} if $V(G)$ and $E(G)$
are finite. In this paper we always assume that all graphs are finite and
have no isolated vertices. 

An edge $e$ is called a {\it loop} if $\partial e=\{v,v\}$ for some vertex $v$.
Two edges $e$ and $e'$ are called {\it parallel} if 
$\partial e=\partial e'=\{a,b\}$, $a\ne b$. A graph is called {\it loop-free}
if it has no loops. A graph is {\it simple} if it is loop-free and
has no parallel edges. A simple graph can be defined as a pair $(V,E)$ where
$E$ is a set of unordered pairs of distinct vertices.

The number of edges incident to a vertex $v$ (loops counted twice)
is called the {\it degree} of $v$.
We say that a graph is {\it of minimum degree} $k$ if the degree of each vertex
is at least $k$. A graph is called {\it $k$-regular} if the degree of every
vertex is $k$.

To each graph we associate
a CW-complex of dimension 1 in the standard way. Usually we shall not 
distinguish between a graph and the corresponding CW-complex.
%
However, when speaking of graphs,
{\it removal of a vertex\/} $v$ always means removal of $v$ together with
all the incident edges. So the result is still a graph with one vertex less
(not the non-compact space obtained by deleting a vertex from the corresponding
CW-complex).

A graph is called $k$-{\it connected\/} (resp. $k$-{\it edge-connected\/})
if the removal of less than $k$ vertices (resp. edges)
cannot disconnect the graph.

\smallskip
\subhead 2.2. Graphs on surfaces
\endsubhead
In this paper, a {\it surface} is a smooth compact 2-manifold without boundary.

Let $G$ be a graph embedded in a surface $F$.
The connected components of $G\setminus F$ are called 
{\it regions}. The pair $(F,G)$ (or just $G$ when it is clear which surface $F$
is considered) is called {\it cellular} if each region is homeomorphic to
an open disk. In this case the regions are called the {\it faces} of $G$.
A pair $(F,G)$ is cellular if and only if $F$ admits a structure of CW-complex
such that $G$ is the 1-skeleton and $V(G)$ is the 0-skeleton.
It is easy to see that any cellular embedded graph is connected and any connected
graph in $S^2$ is cellular.

To avoid any ambiguity between an edge (resp. vertex) and its occurrence
in the boundary of a given region $r$, we call the latter {\it side\/}
(resp. {\it corner\/}) of $r$. In other words, a side (resp. corner) of $r$ is an edge
(resp. vertex) adjacent to $r$ which is considered together with a 
small portion of $r$ near it.
The number of sides of a face $f$ is called the {\it degree} of $f$.
A face of degree $n$ is called also an $n$-gon 
(union, digon, triangle, quadrangle, pentagon, etc. for $n=1,2,3,4,5,\dots$).

A cellular graph is called {\it $2$-cell-embedded} if all sides and corners
of any face are represented by pairwise distinct edges and vertices.

A graph $G$ on a surface $F$ is called {\it simply embedded} if it is loop-free
and for any two parallel edges $\alpha$ and $\beta$, the circle $\alpha\cup\beta$
does not bound a disk in $F$. In particular, a graph in $S^2$ is simply embedded
if and only if it is simple.

If $(F,G)$ is cellular, we define
the {\it dual graph} of $G$ and denote it by $\check G$. It has exactly one vertex
in each face of $G$ and there is a bijection between the edges of $G$ and those
of $\check G$ such that each edge of $\check G$ transversally crosses 
the corresponding edge of $G$ at a single point.
According to the previous definition, the degree of a face of $G$ is equal to the 
degree of the corresponding vertex of $\check G$ and vice versa.

A graph $G$ embedded in a surface $F$ is called a {\it quadrangulation} of $F$
if all its regions are quadrangles.
Note that we do not claim in this definition that $G$ is simple or 2-cell-embedded
(as it is demanded in [\refN] and [\refBGGMTW]).
For example, if $p$ is a point on the circle $S^1$, then
$(S^1\times\{p\})\cup(\{p\}\times S^1)$ is a quadrangulation of the torus
$T=S^1\times S^1$ which has one vertex and two loops.
A 3-path (i.~e. the graph ${{\bullet}\!\!-\!\!-\!\!{\bullet}\!\!-\!\!-\!\!{\bullet}}$)
on a 2-sphere or a non-contractible 2-cycle on $\RP^2$ are also
examples of quadrangulations.

\smallskip
\subhead   2.3.  Basic nets on surfaces
\endsubhead
%
%
A {\it net} on a surface $F$ is the image of a generic immersion of several
circles. In particular, a connected net is either a circle or it can be
represented by a connected 4-regular embedded graph.

\smallskip
{\bf Convention  \remDegFour. }
If $\Gamma$ is a connected net which is not an embedded circle, then
we consider $\Gamma$ as a $4$-regular graph (i.~e., all vertices of $\Gamma$
are crossing points).
\smallskip

A net $\Gamma$ is called {\it irreducible} if for any embedded circle
transversally intersecting $\Gamma$ at most at two points
and dividing $F$ into two components $F_1$ and $F_2$, the following
condition holds: one of $F_1$, $F_2$
is a disk whose intersection with $\Gamma$ is either a simple arc or empty.

A net $\Gamma$ on $F$ is called {\it basic} if it is 
irreducible and none of the components of $F\setminus\Gamma$
is a digon whose corners are at two distinct vertices.

It is easy to check that in the case when $F$ is $S^2$ or
$\RP^2$, our definition of a basic net is equivalent to
the definitions given in [\refC] and [\refD] respectively
(but our definition of an irreducible net differs from that in [\refD]).
Following Conway [\refC], basic nets on $S^2$ are usually called 
{\it basic polyhedra}.

\head 3. Basic properties of basic nets
\endhead
\subhead 3.1. Generalities
\endsubhead

\proclaim{ Proposition \propOne } Let $\Gamma$ be a basic net on a surface $F$. Then:

\smallskip
(a). Any region of $\Gamma$ is planar, i.~e., homeomorphic to a subset of $\R^2$.

\smallskip
(b). If $F$ is a sphere or $\RP^2$, then $\Gamma$ is cellular, in particular,
$\Gamma$ is connected.

\smallskip
(c). If $F=\RP^2$ and $\Gamma$ is an embedded circle,
then $\Gamma$ is a non-contractible curve {\rm(}a pseudoline\/{\rm)}.
\endproclaim

\demo{ Proof } 
(a). Let $r$ be a region of $(F,\Gamma)$. It is an open surface of finite type, thus
$r$ is a connected sum of a planar surface and a compact surface without boundary.
Thus means that there is an embedded circle $\gamma$ which cuts $r$ into two
parts $r_0$ and $r_1$ such that $r_0$ is planar and $\partial r_1=\gamma$.
Since $(F,\Gamma)$ is irreducible, $r_1$ is a disk, hence $r$ is planar.

\smallskip
(b). Let $r$ be a region of $(F,\Gamma)$. By (a), $r$ is planar. Suppose that
$r$ has more than one boundary component. Then $r$ can be cut by an embedded
circle $\gamma$ into two parts such that each part is adjacent to $\Gamma$.
Since $\gamma$ divides $r$, the normal bundle of $\gamma$ is trivial, hence
$\gamma$ divides $F$ which contradicts the irreducibility of $(F,\Gamma)$.

\smallskip
(c). Follows from (b).
\qed
\enddemo


\proclaim{ Proposition \propLoopFree } Let $F$ be either a sphere or a projective plane.
Let $\Gamma$ be a basic net on $F$. Suppose that $\Gamma$ is a $4$-regular graph 
which is not loop-free. Then it has one vertex and
two edges {\rm(}the edges are loops\/{\rm)}.
Moreover, if $F$ is a sphere, then $\Gamma$ is a ``figure-eight''
curve {\rm(}$1^*$ in Figure \figBNS\/{\rm)};
if $F=\RP^2$, then $\Gamma$ is a union of two pseudolines
{\rm(}$g^1$ in Figure \figBNRP\/{\rm)}.
\endproclaim

\demo{ Proof }
Suppose that $\Gamma$ has a loop $\alpha$ adjacent to a vertex $v$.
Let $N$ be a tubular neighbourhood of $\alpha$ in $F$.
It is either an annulus or a M\"obius band.

\smallskip
If $N$ is a M\"obius band, then $\partial N$ is an embedded circle
intersecting $\Gamma$ at two points. Since $\Gamma$ is irreducible,
$\partial N$ bounds a disk whose intersection with $\Gamma$ is a simple arc.
Then the edge of $\Gamma$ containing this arc is another loop $\beta$ adjacent to $v$.
and the result follows (if $\beta$ were not a pseudoline, then $\Gamma$ would
be reducible).

\smallskip
Now suppose that $N$ is an annulus. Since $F$ is
a sphere or a projective plane, each of the two components of $\partial N$
divides $F$ and intersects $\Gamma$ at most at two points.
Hence, the irreducibility of $\Gamma$ implies that $F\setminus N$ is a union of
two disjoint disks and $\Gamma\setminus N$ is a simple arc on one of them.
\qed\enddemo

\proclaim{ Lemma \lemZero }
Let $\Gamma$ be a net in a surface $F$ and let $\alpha$ be a simple closed
curve transversal to $\Gamma$ and null-homologous in $F$.
Then $\alpha$ cuts $\Gamma$ at an even number of points. \qed
\endproclaim

\subhead 3.2. Basic nets in $S^2$
\endsubhead
The following fact is well-known but we give a precise
statement and a proof for the sake of completeness.

\proclaim{ Proposition \propS } Let $\Gamma$ be a connected graph on $S^2$ and $G$ its dual.
The following conditions are equivalent:
\roster
\item
 $\Gamma$ is a basic net which is neither a circle nor a figure-eight curve; 
\item
 $\Gamma$ is simple, $4$-regular, and $4$-edge-connected;
\item
 $G$ is a simple quadrangulation of minimum degree $3$;
\item
 $G$ is a simple, $2$-cell-embedded, $2$-connected, and $3$-edge-connected
 quadrangulation of minimum degree $3$.
\endroster
\endproclaim

\demo{ Proof } Note that $\Gamma$ is $4$-regular if and only if $G$ is a quadrangulation.
So, we assume from now on that $\Gamma$ is $4$-regular and $G$ is a quadrangulation.
\smallskip

(1) $\Longrightarrow$ (2).
Suppose that Condition (1) holds.

{\it Simplicity.}
$\Gamma$ is loop-free by Proposition \propLoopFree.
Let us prove that $\Gamma$ cannot have parallel edges. Suppose that
$\alpha$ and $\beta$ are two parallel edges. 
Then $\alpha\cup\beta$ is a simple closed curve (not necessarily smooth).
Let $N$ be a tubular neighbourhood of $\alpha\cup\beta$.
It is an annulus and the pair $(N,N\cap\Gamma)$ is as in one of
Figures \figAnnuli.1(a--d).
Case (d) is impossible by Lemma \lemZero.
In Cases (a--c), the irreducibility condition implies that one
of the components of $\partial N$ (the interior one in Figure
\figAnnuli.1) bounds a disk $D$ such that $D\cap\Gamma$ is $\varnothing$
or an arc. In Cases (a) and (c) this provides a digon and
in Case (b) this contradicts the irreducibility (see Figure \figAnnuli.2).
So, we proved that the graph $\Gamma$ is simple.

\midinsert
\centerline{
  \epsfxsize=16mm
  (a)\!\!\!\epsfbox{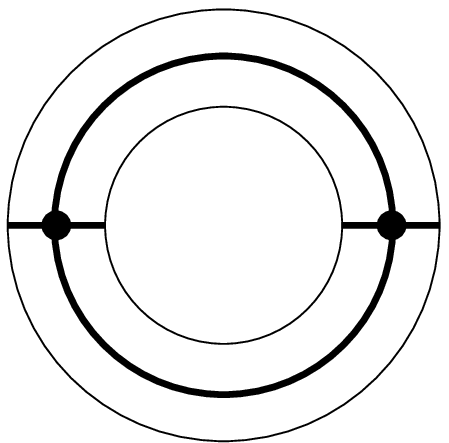}\hskip2mm
  \epsfxsize=16mm
  (b)\!\!\!\epsfbox{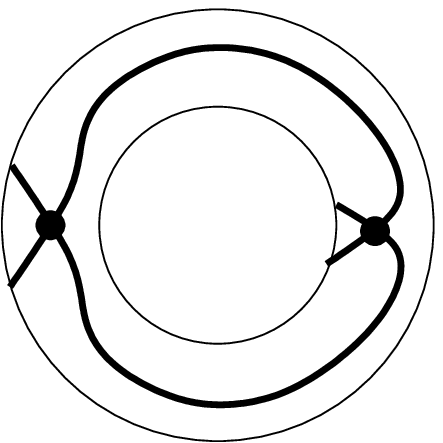}\hskip2mm
  \epsfxsize=16mm
  (c)\!\!\!\epsfbox{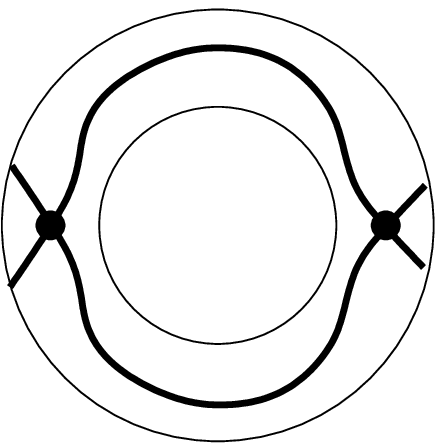}\hskip2mm
  \epsfxsize=16mm
  (d)\!\!\!\epsfbox{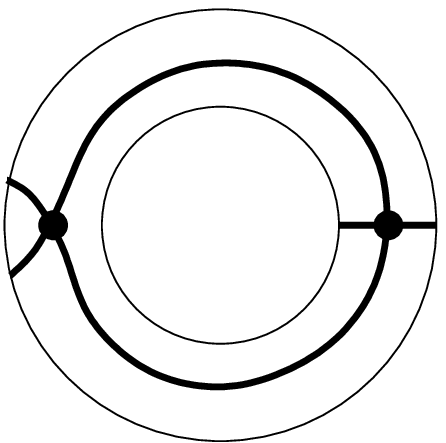}
\hskip12mm
  \epsfxsize=16mm\epsfbox{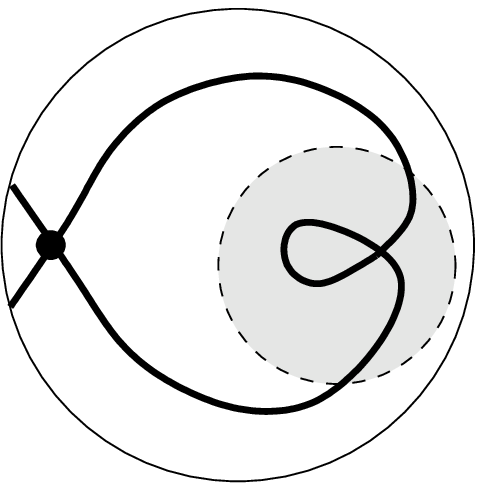}
}
\botcaption{
\hskip25mm
Figure \figAnnuli.1. \hskip42mm
Figure \figAnnuli.2.}
\endcaption
\endinsert

{\it $4$-edge-connectivity.}
Let $e_1,\dots,e_k$ be a minimal set of edges which disconnects $\Gamma$.
Then $\Gamma\setminus\bigcup_i e_i$ has two connected components
$\Gamma_1$ and $\Gamma_2$,
and each edge $e_i$ relates them. Hence there exists an embedded circle
$\gamma$ which separates $\Gamma_1$ from $\Gamma_2$ and
transversally crosses every edge $e_i$ at one point.
Hence $k$ is even by Lemma \lemZero. Since $\Gamma$ is irreducible,
$k$ cannot be $0$ or $2$. Thus, $k\ge4$, i.~e., $\Gamma$ is $4$-edge-connected.

\smallskip
(2) $\Longrightarrow$ (1).
Indeed,
the $4$-edge-connectivity easily implies the irreducibility
and the simplicity (the absence of parallel edges) implies
the absence of digons.

\smallskip
(2) $\Longrightarrow$ (3).
Suppose that $\Gamma$ is simple and $4$-edge-connected.

%
A loop of $G$ would cut $\Gamma$ at one point which is impossible by Lemma \lemZero,
hence $G$ is loop-free. Suppose that $\alpha$ and $\beta$ are parallel edges of $G$.
Then $\alpha\cup\beta$ is a circle which cuts $\Gamma$ at two points. Since
$\Gamma$ is irreducible, these two points are connected by a simple arc of $\Gamma$,
hence they belong to the same edge of $\Gamma$ which contradicts the definition of
the dual graph. Thus, $G$ is simple.

Let us show that the minimum degree is $3$.
Indeed, let $v\in V(G)$. If $\deg(v)=1$, then the edge adjacent to $v$ is dual to a loop of $\Gamma$.
If $\deg(v)=2$, then the face dual to $v$ is a digon.

\smallskip
(3) $\Longrightarrow$ (4).
Suppose that $G$ is simple of minimum degree $3$.
Let us show that it is {\it $2$-cell-embedded}.
Indeed, let $f$ be a face of $G$.
Suppose that two sides of $f$ are represented by the same edge $e$.
If they are consecutive, i.~e., if they have a common corner at $v$, then $\deg(v)=1$.
If they are opposite, then each of the two other sides represents a loop.
So, we conclude that all sides of $G$ are represented by pairwise distinct edges.
Suppose that $f$ has two corners at the same vertex.
If they are consecutive, then the side between them is a loop (see Figure \figPropS(a)).
If they are opposite, then $G$ has parallel edges (see Figure \figPropS(b)).

\midinsert
\centerline{
  \epsfxsize=18mm
  (a)\!\!\!\epsfbox{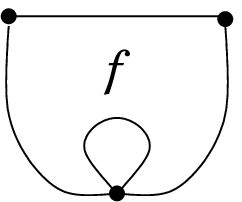}\hskip10mm
  \epsfxsize=16mm
  (b)\!\epsfbox{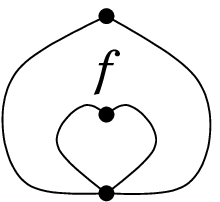}\hskip2mm
}
\botcaption{ Figure \figPropS }
\endcaption
\endinsert

$G$ is {\it $2$-connected\/} because otherwise it would not be $2$-cell-embedded.

Let us show that $G$ is {\it $3$-edge-connected}. Indeed, let
$e_1,\dots,e_k$ be a minimal set of edges which disconnects $G$. Then  there exists an embedded
circle $\gamma$ which transversally crosses every edge $e_i$ at one point.
We have $k\ge2$ because $G$ is $2$-cell-embedded.
Suppose that $k=2$. Let $f$ and $f'$ be the faces of $G$ crossed by $\gamma$.
If $e_1$ and $e_2$
have a common vertex $v$, then (since $\deg v>2$) $v$ represents two corners of one of
the faces $f$ or $f'$.
So, we conclude that $e_1$ and $e_2$ do not have a common vertex. Then $e_1$ and $e_2$ are
opposite sides of both faces $f$ and $f'$.
Let $u$ and $v$ be the ends of $e_1$ and $e_2$ on the
same side of $\gamma$. Then $G$ has an edge $uv$ which is a common side of $f$ and $f'$.
Hence $\deg u=\deg v=2$. Contradiction.

\smallskip
(4) $\Longrightarrow$ (2).
Suppose that $G$ is simple, $3$-edge-connected of minimum degree $3$.
If $\Gamma$ has a loop, then the removal of its dual edge disconnects $G$.
If $\Gamma$ has two parallel edges, then the removal of their duals disconnects $G$.
Thus, $\Gamma$ is simple.
The edge-connectivity of $\Gamma$ is even by Lemma \lemZero\
and it cannot be equal to 2 because $G$ has no parallel edges.
\qed
\enddemo


\subhead  3.3. Basic nets in $\RP^2$
\endsubhead

\proclaim{ Theorem \thRP }
Let $\Gamma$ be a cellular graph in $\RP^2$ and $G$ its dual.
Let $\xi:S^2\to\RP^2$ be the universal covering and let
$\tilde\Gamma = \xi^{-1}(\Gamma)$. 
Then the following conditions are equivalent:
\roster
\item
  $\Gamma$ is a basic net which is neither a line nor a union of two lines;
\item
  $\Gamma$ is simply embedded, $4$-regular, $4$-edge-connected, and $G$ is loop-free;
\item
  $G$ is a simply embedded quadrangulation of minimum degree $3$;
\item
  $\tilde\Gamma$ is a basic net in $S^2$ which is not a circle.
\endroster
\endproclaim

\demo{ Proof }
(1) $\Longrightarrow$ (2). Assume that Condition (1) holds.
Let us show that $G$ is loop-free.
Indeed, if $\alpha$ is a loop of $G$, then it cuts $\Gamma$ at one point.
By Lemma \lemZero, this implies that $\alpha$ is a pseudoline. Let $N$ be a tubular
neighbourhood of $\alpha$ and let $D=\RP^2\setminus N$.
Then $N$ is a M\"obius band and $\Gamma\cap N$ is a simple arc, in particular,
$\partial N$ is a circle which cuts $\Gamma$ at two points.
Then $D\cap\Gamma$ is a also simple arc because $\Gamma$ is irreducible, i.~e.,
$\Gamma$ is a circle which is impossible by Condition (1).
Thus, $G$ is loop-free.
The rest of the proof is the same as in Proposition \propS.

\smallskip

\if01{
(2) $\Longrightarrow$ (1). Assume that Condition (2) holds. Then:

\smallskip
{\it $\Gamma$ is irreducible.}
Indeed, let $\gamma$ be an embedded circle which divides $\RP^2$.
Then $\RP^2\setminus\gamma=M\cup D$ where $D$ is a disk and $M$ is a M\"obius band.
Suppose that $\gamma$ transversally meets $\Gamma$ at $k\le 2$ points.
Note that $k$ is even by Lemma \lemZero.
If $k=0$, then (since $\Gamma$ is connected) either $\Gamma\cap D$ or $\Gamma\cap M$ is empty,
but $\Gamma\cap M=\varnothing$ is impossible because $\Gamma$ is cellular.
If $k=2$, then the both points of $\Gamma\cap\gamma$ must belong to the same edge
because $\Gamma$ is $4$-edge-connected. Then either $D\cap\Gamma$ or $M\cap\Gamma$
is simple arc, but the later case is impossible because $G$ is loop-free.

\smallskip
{\it No face of $\Gamma$ is a digon} because $\Gamma$ is simply embedded.

\smallskip
}\fi

(2) $\Longrightarrow$ (3). The same proof as in Proposition \propS\
(note that $G$ is already loop-free by Condition (2)).

\smallskip

(3) $\Longrightarrow$ (4). It is immediate to check that
if $G$ is a simply embedded quadrangulation of $\RP^2$
of minimum degree 3, then $\xi^{-1}(G)$ is a simple quadrangulation of $S^2$
of minimum degree 3. Thus, the result follows from Proposition \propS.

\smallskip

(4) $\Longrightarrow$ (1).
Assume that $\tilde\Gamma$ is a basic net in $S^2$ 
and let us prove that $\Gamma$ is a basic net in $\RP^2$.
If $D$ is a digon of $\Gamma$, then $\xi^{-1}(D)$ is a digon of $\tilde\Gamma$,
thus it remains to prove that $\Gamma$ is irreducible. Indeed,
let $\gamma$ be an embedded circle transversally intersecting $\Gamma$
at $i\le2$ points and dividing $\RP^2$ into two components
($i=0$ or $2$ by Lemma \lemZero).
One of the components 
is an open disk $D$. Let $\xi^{-1}(D)=\tilde D_1\sqcup\tilde D_2$
and $\tilde\gamma_1=\partial\tilde D_1$.
If $i=0$, then $D\cap\Gamma=\varnothing$ because $\tilde\Gamma$ is connected.
Let $i=2$. 
Since $\tilde\Gamma$ is irreducible, $\tilde D\cap\tilde\Gamma$ is a simple arc
where $\tilde D$ is one of the two components of $S^2\setminus\tilde\gamma_1$.
If $\tilde D=\tilde D_1$, then $D\cap\Gamma=\xi(\tilde D\cap\tilde\Gamma)$ is a simple arc
and we are done. Otherwise $\tilde\Gamma\setminus\tilde D_1$ is a simple arc, hence
its subset $\tilde\Gamma\cap\tilde D_2$ is a priori a disjoint union of simple arcs,
but the total number of their boundary points is $2$, hence it
is a simple arc.
\qed
\enddemo

\head  \sectCover. Uniqueness of a planar projective quotient of a planar graph.
\endhead

\proclaim{ Theorem \thUniq }
Let $G_1$ and $G_2$ be embedded graphs in $\RP^2$ without vertices of degree $2$.
Let $\xi_1$ and $\xi_2$ be two unramified coverings $S^2\to\RP^2$ and
let $\sigma_j:S^2\to S^2$, $j=1,2$, be the corresponding deck transformations,
i.~e., for any $x\in S^2$, $\sigma_j(x)=y$ where $y\ne x$ and $\xi_j(x)=\xi_j(y)$.
Suppose that $\xi_1^{-1}(G_1) = \xi_2^{-1}(G_2)$ and that it is a connected graph
{\rm(}we denote it by $G$\/{\rm)}. Then
$\sigma_1|_G$ and $\sigma_2|_G$ are combinatorially equivalent, i.~e.,
$\sigma_1|_{V(G)} = \sigma_2|_{V(G)}$ and
for any $e\in E(G)$ we have $\sigma_1(e)=\sigma_2(e)$.
\endproclaim

\demo{ Proof }
Without loss of generality we may assume that $S^2$ is glued out of regular
polygons ($G$ being represented by their sides) and the mappings $\sigma_j$ are
linear on each of them. Then $\sigma_1$ and $\sigma_2$ are combinatorially equivalent
if and only if $\sigma_1=\sigma_2$.

We set $\tau=\sigma_1\circ\sigma_2$.
Since $\sigma_1^2=\sigma_2^2=\id_{S^2}$, it is enough to prove that $\tau=\id_{S^2}$.
Note that $\tau$ is an orientation preserving homeomorphism $S^2\to S^2$.
We suppose that $\tau\ne\id$ and we shall obtain a contradiction in several steps.

\smallskip
{\bf Step 1.} {\it There do not exist $v\in V(G)$ and an edge
$e$ adjacent to $v$ such that
$\tau(v)=v$ and $\tau(e)=e$.}
Indeed, $\tau$ is the identity map on the faces adjacent to $v$.
The same is true for faces adjacent to them etc. Since the graph $G$ is
connected, we exhaust all its vertices and edges by this process.

\smallskip
{\bf Step 2.} {\it $\tau$ has exactly two fix points.}
Indeed, let $L(\tau)$ be the Lefschetz number of $\tau$, i.~e.,
$L(\tau)=\sum_q (-1)^q\trace(\tau_*:H_q(S^2)\to H_q(S^2))$.
Since $\tau$ is an orientation preserving homeomorphism, we have $L(\tau)=2$.
It is well-known that $L(\tau)$ is equal to the intersection number of
the diagonal of $S^2\times S^2$ with the graph of $\tau$.
We deduce from the result of Step 1 that the number of fix points of $\tau$ is finite.
Moreover, at any fix point, $\tau$ is locally conjugated to a rotation, hence
the local intersection of the diagonal with the graph at any fix point is equal to $+1$.

\smallskip
{\bf Step 3.} {\it We denote the fix points of $\tau$ by $x$ and $y$.
Then $\sigma_1(x)=\sigma_2(x)=y$ and $\sigma_1(y)=\sigma_2(y)=x$.}
Indeed, let $z=\sigma_2(x)$. By the definition of $\tau$ we have
$\tau(z)=\sigma_1(\sigma_2(\sigma_2(x)))=\sigma_1(x)$.
Since $x$ is a fix point of $\tau$, we have
$\sigma_1(x)=\sigma_1(\tau(x))=\sigma_1(\sigma_1(\sigma_2(x)))=\sigma_2(x)=z$.
Thus $\tau(z)=\sigma_1(x)=z$, i.~e., $z$ is a fix point of $\tau$. Hence $z=y$.
The other equalities are obtained similarly.

\smallskip
{\bf Step 4.}
Subdividing if necessary the faces containing $x$ and $y$ we may assume
without loss of generality that $x$ and $y$ are vertices of $G$.
Let $\gamma$ be a shortest path on $G$ from $x$ to $y$
(a path with the minimum number of edges). Let
$x=x_0, x_1,\dots,x_n=y$ be the successive vertices on $\gamma$.
Then $\sigma_1(\gamma)\cap\sigma_2(\gamma)=\{x,y\}$. Indeed,
Suppose that $\sigma_1(x_i)=\sigma_2(x_j)$. If $i\ne j$, say, $i<j$,
then $\gamma$ is not a shortest path from $x$ to $y$ because
in this case the path $x=\tau(x_0),\tau(x_1),\dots,\tau(x_i)=x_j,
x_{j+1},\dots,x_n=y$ is yet shorter. Hence $i=j$ and so $x_i$ is a fix point
of $\tau$, hence $x_i=x_j\in\{x,y\}$.

Thus, $\sigma_2(\gamma)$ is contained in one of the two disks bounded by the
circle $\gamma\cup\sigma_1(\gamma)$. This contradicts the fact that
each of the circles  $\gamma\cup\sigma_1(\gamma)$ and $\gamma\cup\sigma_2(\gamma)$
divides the sphere into two halves containing the same number of 2-faces.
\qed\enddemo

\head \sectGener. Generating basic nets
\endhead

In this section we prove Theorem 2.
By the duality (see Theorem \thRP), it follows easily from Theorem \thMin\
(see below) combined with Theorem 1.

\subhead \sectGener.1. Face contraction/removal
\endsubhead
We recall here some definitions from [\refN], [\refBGGMTW].
\smallskip
Let $G$ be a quadrangulation of a surface $F$ which has more than one face
and let $f=abcd$ be a face of $G$ such that $a\ne c$. In this case we say that
the face $f$ is {\it contractible at} $\{a,c\}$ and the {\it contraction of $f$
at} $\{a,c\}$ consists in the removal of the interior of $f$ and glueing the edges
$ba$ with $bc$ and $da$ with $dc$ (see Figure \figContract).
The inverse operation to a face contraction is called a {\it vertex splitting}.

\midinsert
\epsfxsize=100mm
\centerline{\epsfbox{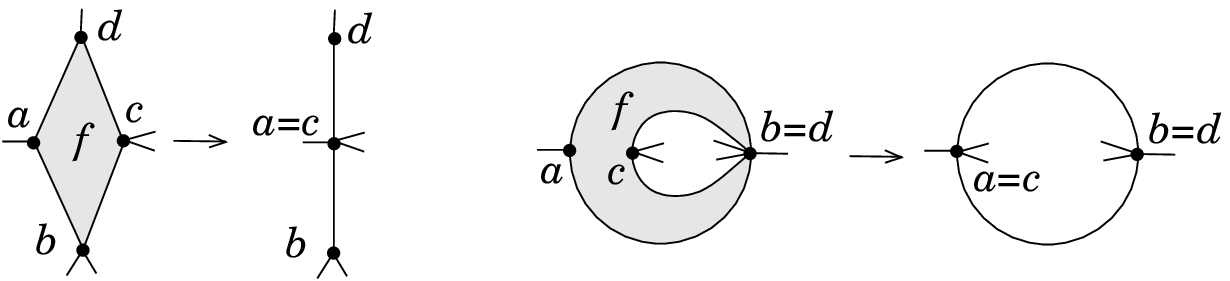}}
\botcaption{ Figure \figContract}
   Examples of face contraction
\endcaption
\endinsert

\smallskip
Let $G$ be a quadrangulation of a surface $F$ 
and let $f=abcd$ be a face of $G$. We say that $f$ is {\it removable} if
$a,b,c,d$ are pairwise distinct vertices of degree $3$ and,
if we denote their outcoming edges (not being the sides of $f$) by $aa_1,bb_1,cc_1,dd_1$,
then $\{a,b,c,d\}\cap\{a_1,b_1,c_1,d_1\}=\varnothing$.
In this case, the
{\it removal of} $f$ consists just in the removal of the vertices $a,b,c,d$
and all the edges incident to them (see Figure \figRemoval).
The inverse operation to a face removal is called a {\it face addition}.

\midinsert
\epsfxsize=120mm
\centerline{\epsfbox{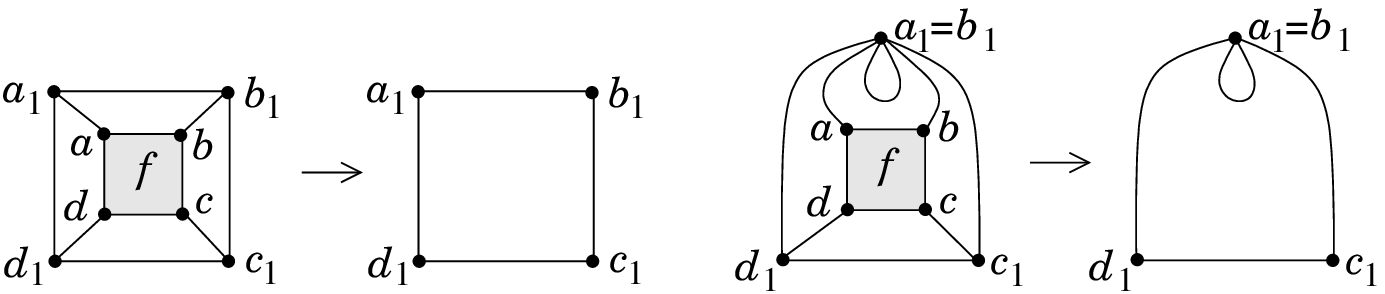}}
\botcaption{ Figure \figRemoval}
   Examples of face removal
\endcaption
\endinsert

It is easy to see that the result of a face contraction/removal is again a
quadrangulation of the same surface.

\smallskip
Let $\Cal Q$ (resp. $\bar\Cal Q$) be the class of simple (resp. simply embedded)
quadrangulations of $S^2$ (resp. of $\RP^2$) of minimal degree $3$.
By Proposition \propS\ (resp. by Theorem \thRP), it is
the dual of the class of basic nets on $S^2$ (resp. on $\RP^2$) with more
than one crossing.

\smallskip
Given $G\in\Cal Q$, we say that a face $abcd$ of $G$ is $\Cal Q$-removable
(resp. $\Cal Q$-contractible
at $\{a,c\}$) if it is removable (resp. contractible at $\{a,c\}$) and the
result of the removal (resp. contraction) belongs to $\Cal Q$.
We say that a $\Cal Q$-contraction of a face $abcd$ at $\{a,c\}$
is {\it special\/} if $\deg a=3$ or $\deg c=3$.

If $G\in\Cal Q$ does not have any special
$\Cal Q$-contractible or $\Cal Q$-removable face, then
we say that $G$ is {\it $\Cal Q$-minimal}.

In the same way we define (special) $\bar\Cal Q$-contractible/removable faces
(of quadrangulations belonging to $\bar\Cal Q$)
and $\bar\Cal Q$-minimal quadrangulations. Let 
$$
   \Cal Q_{\min}=\{G\in\Cal Q\mid\text{$G$ is $\Cal Q$-minimal}\},\quad
   \bar\Cal Q_{\min}=\{G\in\bar\Cal Q\mid\text{$G$ is $\bar\Cal Q$-minimal}\},
$$

It is clear that (special) $\Cal Q$- or $\bar\Cal Q$-vertex-splittings and
$\Cal Q$- or $\bar\Cal Q$-face-additions are dual to (special) 
face splittings and vertex surroundings on basic nets
respectively.

\subhead\sectGener.2. Double covering and minimality
\endsubhead
Let $\xi:S^2\to \RP^2$ be the double covering and
$\sigma:S^2\to S^2$ its deck transformation,
i.~e., $\xi\circ\sigma=\xi$ and $\sigma\ne\id$.
For $a\in S^2$ or $a\subset S^2$,
we denote $\sigma(a)$ by $a'$ and $\xi(a)$ by $\bar a$.

\proclaim{ Lemma \lemMin } Let $\bar G\in\bar\Cal Q$
and $G=\xi^{-1}(\bar G)$. Suppose that a face $f=abcd$ of
$G$ is $\Cal Q$-contractible at $\{a,c\}$. Then:


(a) The vertices $a$, $a'$, $c$, $c'$ are pairwise distinct.

(b) $\bar f$ is not $\bar\Cal Q$-contractible at $\{\bar a,\bar c\}$
if and only if one of the following two conditions holds:
{\rm(i)} $G$ has an edge $ac'$ or
{\rm(ii)} $b=d'$ and $\deg b=4$.

(c) Assume, moreover, that $f$ is special $\Cal Q$-contractible at $\{a,c\}$.
Then $\bar f$ is special $\bar\Cal Q$-contractible at $\{\bar a,\bar c\}$
if and only if $ac'\not\in E(G)$.
\endproclaim

\demo{ Proof }
(a). We have $a\ne c$ because $G\in\Cal Q$ by Theorem \thRP,
hence $G$ is 2-cell-embedded by Proposition \propS.
We have $a\ne c'$ because
otherwise we have $f'=a'b'c'd'=cb'ad'$ and the result of the contraction of $f$ is 
not 2-cell-embedded at the face
$f'$ which contradicts Proposition \propS.

\smallskip
(b). By (a) we have $\bar a\ne\bar c$, i.~e., the face $\bar f$ is contractible at
$\{\bar a,\bar c\}$. Let $\bar G_1$ be the result of the contraction.

If (i) holds, then the image of the edge $\bar a\bar c$ on $\bar G_1$ would be a loop;
if (ii) holds, then $\deg_{\bar G_1}\bar b=2$.
In both cases we have $\bar G_1\not\in\bar\Cal Q$.

Suppose that none of Conditions (i), (ii) holds.
Let us show that $\bar G_1\in\bar\Cal Q$. It is clear that $\bar G_1$ is a quadrangulation.
Since $f$ is $\Cal Q$-contractible, we have $\deg\bar b=\deg b>3$ and $\deg\bar d=\deg d>3$.
Since Condition (ii) does not hold, it follows that
$\bar G_1$ is a quadrangulation of minimum degree 3.
So, it remains to prove that $\bar G_1$ is simply embedded.
Since $\bar G$ is loop-free and $a\ne c'$, it follows that $\bar G_1$ is loop-free also.
Suppose that there are two parallel edges on $\bar G_1$ which bound a disk $D^*$ on $\RP^2$.
Let $\pi:\RP^2\to\RP^2$ be a continuous mapping which extends the contraction of $\bar f$ so that
$\pi|_{\bar f}$ is constant on each segment parallel to the diagonal $\bar a\bar c$ and
$\pi|_{\RP^2\setminus\bar f}$ is a homeomorphism onto its image.
Then $\bar D=\pi^{-1}(D^*)$ is a disk bounded by two edges of $\bar G$ and, maybe, by the diagonal
of $\bar f$ if $\pi(\bar a)=\pi(\bar c)\in\partial D^*$.
Thus, $\xi^{-1}(\bar D)$ is a disjoint union of two disks $D\cup D'$ on $S^2$ such that
$D$ is bounded by two edges of $G$ and, maybe, by the diagonal $ac$ of the face $f$.
Thus either $G$ or the result of the contraction of $f$ at $\{a,c\}$ is not simple. Contradiction.

\smallskip
(c). Since $f$ is special $\Cal Q$-contractible at $\{a,c\}$,
without loss of generality we may assume that
$\deg(a)=3$. By (b), it is enough to show that Condition (ii) does not hold. Suppose that
it does hold. Then, by (a), $G$ contains a subgraph depicted in Figure \figB.1.
Here we denote the third outcoming edge from $a$ by $ax$. It is clear that $x$
should be in the quadrangle $q=abc'd$. It cannot be on the boundary of $q$. Indeed,
$x\not\in\{a,b,d\}$ because $G$ is simple and $x\ne c'$ because $G$ is bipartite (see
the colors in Figure \figB.1). We have $\deg(b)=\deg(d)=4$ and $\deg(a)=3$, hence
all outcoming edges from $a$, $b$, and $d$ are already present in Figure \figB.1.
Therefore, the path $xabc'$ follows the boundary of the same face (we denote it by $f_1$).
Hence $G$ has an edge $xc'$ adjacent to $f_1$. Similarly, $xc'$ is adjacent to the
face $f_2=xadc'$. Since $f_1\cup f_2=q$, we conclude that $\deg(x)=2$. Contradiction.
\qed\enddemo

\midinsert
\epsfxsize=100mm
\centerline{\epsfbox{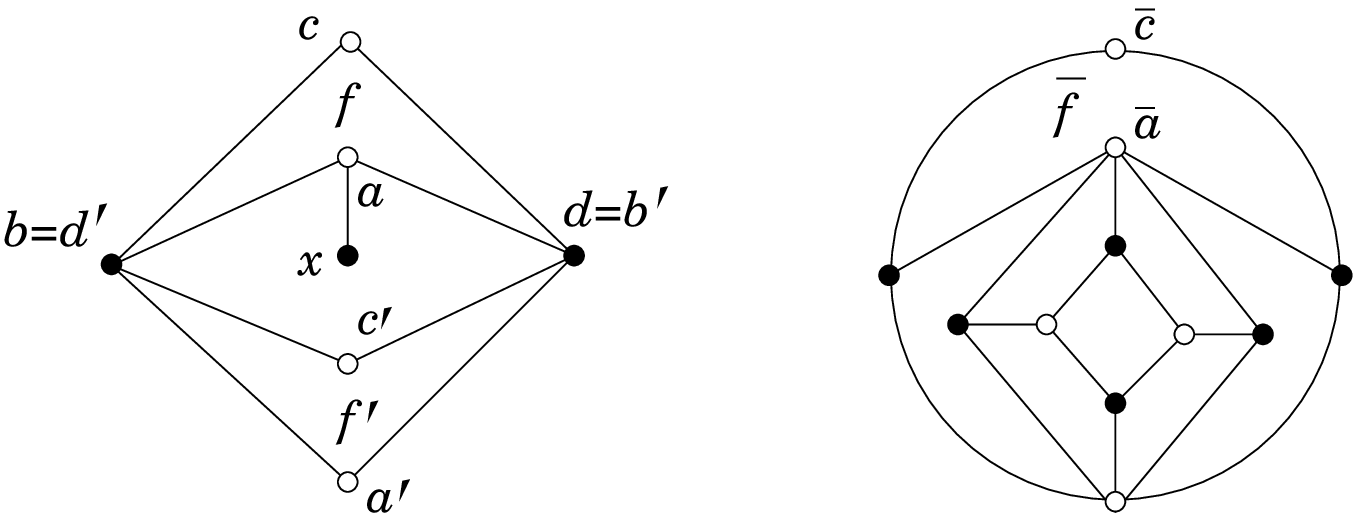}}
\botcaption{\hskip 10mm
Figure \figB.1 \hskip33mm
Figure \figB.2  $\check g^8_5$
}
\endcaption
\endinsert

\definition{ Remark } Condition (ii) of Lemma \lemMin(b) holds for
the quadrangulation $\bar G=\check g_5^8$ depicted in Figure \figB.2 which is dual to 
the basic net $g_5^8$ in Figure \figBNRP. In this case $f$ is (non-special)
$\Cal Q$-contractible at $\{a,c\}$ but $\bar f$ is not $\bar\Cal Q$-contractible at
$\{\bar a,\bar c\}$.
\enddefinition

\proclaim{ Lemma \lemRemove } Let $\bar G\in\bar\Cal Q$
and let $G=\xi^{-1}(\bar G)$. Suppose that a face $f_0$ of
$G$ is removable. Let $f_1,f_2,f_3,f_4$ be the faces which
have a common edge with $f_0$. If $f_i\ne f'_j$ for any $i,j\in\{0,\dots,4\}$,
then $\bar f_0$ is removable.
\endproclaim

\demo{ Proof } Suppose that $\bar f_0$ is not removable.
Than $\bar x=\bar y$ for a vertex $x$ of $f_0$ and
for a vertex $y\ne x$ of one of $f_0,\dots,f_4$.
Since $x\ne y$ and $\bar x=\bar y$, it follows that $y'=x$.
One of the faces adjacent to $y$ is $f_i$ for some $i=0,\dots,4$. Then $f'_i$ is
adjacent to $y'=x$. Since any face adjacent to $x$ is one $f_0,\dots,f_4$, it follows
that $f_i'=f_j$ for some $j=0,\dots,4$.
\qed\enddemo

\midinsert
\epsfxsize=25mm
\centerline{\epsfbox{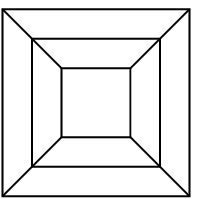}}
\botcaption{Figure \figGfive} The covering of $\check g_2^5$
\endcaption
\endinsert

\proclaim{ Theorem \thMin }
Let $\bar G\in\bar\Cal Q_{\min}$
and $G=\xi^{-1}(\bar G)$.  Then either $G\in\Cal Q_{\min}$ or
$G$ is as in Figure \figGfive\ {\rm(}and then $\bar G=\check g_2^5$ --
the dual of the basic net $g_2^5$ in
Figure \figBNRP\/{\rm)}.
\endproclaim

\demo{ Proof }
Suppose that $G\not\in\Cal Q_{\min}$.
Then it admits either a special $\Cal Q$-face-contraction or a $\Cal Q$-face-removal.

\smallskip
Case 1. $G$ admits a special $\Cal Q$-contraction of a face $f=abcd$ at $\{a,c\}$.
We may assume that $\deg(a)=3$. Since $\bar G\in\bar\Cal Q_{\min}$, the face
$\bar f$ is not $\bar\Cal Q$-contractible at $\{\bar a,\bar c\}$. By Lemma \lemMin(c),
this implies that $G$ has edges $ac'$ and $ca'$ (note that $a\ne c'$ by Lemma \lemMin(a)).
Then $b'\ne d$ because $G$ is bipartite, hence $G$ contains a subgraph $H$
shown in Figure \figA.1.

\midinsert
\epsfxsize=100mm
\centerline{\epsfbox{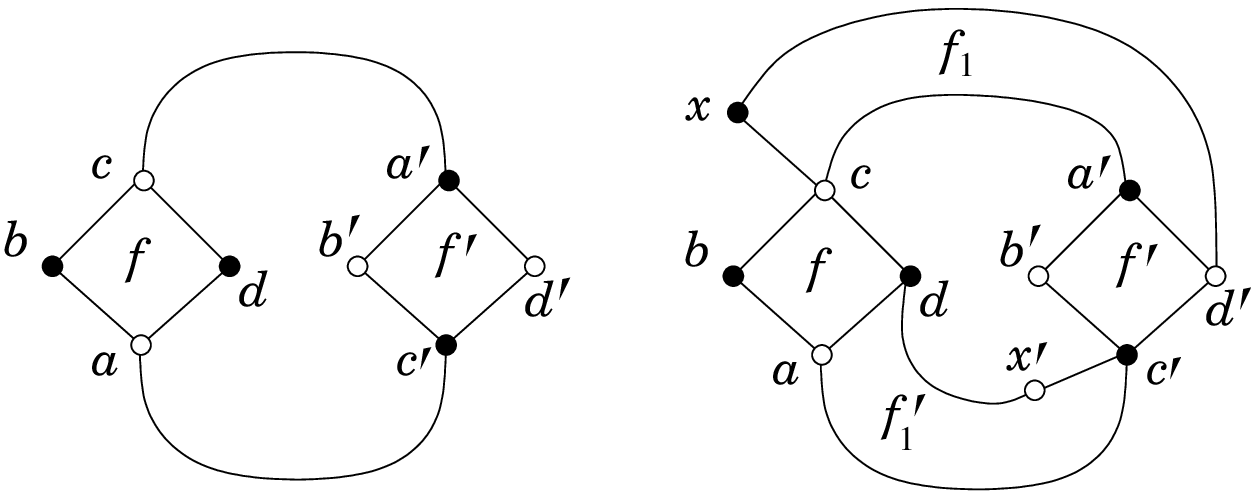}}
\botcaption{
Figure \figA.1 \hskip30mm
Figure \figA.2 
}
\endcaption
\endinsert

Subcase 1.1. $\deg(c)=3$. In this case the paths
$bac'd'$, $dac'b'$, $dca'b'$, and $bca'd'$ belong to the boundaries of some faces.
Hence there are edges $bd'$ and $b'd$ and when we add them to $H$, we
complete the graph $G$. Then $\deg(b)=\deg(d)=3$ which contradicts the condition
that $f$ is $\Cal Q$-contractible at $\{a,c\}$.

Subcase 1.2. $\deg(c)>3$. Let $e=cx$ be an edge adjacent to $c$ which is not in
$H$ and which is next to $ca'$ in the natural cyclic order on the set of
edges outcoming from $c$. Since $H$ is symmetric, we may assume that $e$
sits in the hexagon $h=abca'd'c'$ (the exterior region in Figure \figA.1).
The vertex $x$ is not on the boundary of $h$. Indeed, $x\not\in\{a,d'\}$ because
$G$ is bipartite (see the colors in Figure \figA.1), $x\not\in\{b,a'\}$ because
$G$ is simple, and $x\ne c'$ because $\bar G$ is loop-free.
Thus, $xca'd'$ is a path in $G$. Moreover, by the assumption that
$\deg(a)=3$, this path belongs to the boundary of some face $f_1$ (see Figure
\figA.2).
Then $\bar f_1$ is $\bar\Cal Q$-contractible at $\{\bar a,\bar x\}$ by Lemma \lemMin.
Indeed, we have $\deg(c)>3$ by hypothesis, $\deg(d')=\deg(d)>3$ because $f$ is
$\Cal Q$-contractible at $\{a,c\}$, and there is no edge $xa$ because $\deg(a)=3$
and we have already three outcoming edges from $a$.
Moreover, $\bar f_1$ is special $\bar\Cal Q$-contractible at $\{\bar a,\bar x\}$
because $\deg\bar a=\deg a=3$.

\smallskip
Case 2. $G$ has a $\Cal Q$-removable face $f=abcd$. Let $a_1,b_1,c_1,d_1$ be as
in the definition of the face removal (see \S\sectGener.1) and let $G_1$ be the
result of the removal of the face $f$. Since $G_1$ is in $\Cal Q$, it is
$2$-cell-embedded by Proposition \propS, hence $a_1,b_1,c_1,d_1$ are
pairwise distinct. 
Since $f$ is $\Cal Q$-removable, we have
$$
    \deg x>3 \qquad\text{for $x\in\{a_1,b_1,c_1,d_1\}$}.
                                            \eqno(\eqDeg)
$$

Let us prove that $\bar f$ is removable.
By Lemma \lemRemove\ and by
symmetry, it suffices to check that $a'\not\in\{a,b,c,d,a_1,b_1,c_1,d_1\}$.
We have $a'\not\in\{a_1,b_1,c_1,d_1\}$ by (\eqDeg), $a'\ne a$ because $\sigma$
has no fix point, and $a'\ne b$ because $\bar G$ is loop-free. Suppose that
$a'=c$. Then $b'$ is connected to $c$ by an edge. i.~e., $b'\in\{b,d,c_1\}$.
We have $b'\ne b$ (no fix point of $\sigma$) and $b\ne c_1$ by (\eqDeg), hence
$b'=d$. Thus, $\sigma$ maps the edge $ab$ to the edge $cd$. If follows that
the face $f'$ is incident to $cd$. This is impossible because $f'\ne f$ (otherwise
$\sigma$ has a fix point) and $f'\ne cdd_1c_1$ by (\eqDeg).
So, we proved that $\bar f$ is removable.

Let $\bar G_1$ be the result of
the removal of the face $\bar f$.
Then $\bar G_1$ is a simply embedded quadrangulation
and $\deg_{\bar G_1}(x)\ge 3$ for $x\not\in\{\bar a_1,\bar b_1,\bar c_1,\bar d_1\}$.
Since $\bar G\in\bar\Cal Q_{\min}$, we know that $\bar G_1\not\in\bar\Cal Q$.
Hence the degree in $\bar G_1$
of one of $\bar a_1,\bar b_1,\bar c_1,\bar d_1$, (say, $\bar a_1$) is less that $3$.
Since $\deg_{\bar G}(\bar a_1)=\deg_G(a_1)>3$ (see (\eqDeg)), this means that $\bar a_1$ is
incident in $\bar G$ to at least two edges which are removed in $\bar G_1$.
This may happen only if $a'_1\in\{a_1,b_1,c_1,d_1\}$.
We have
$a'_1\ne a_1$ (since $\sigma$ has no fix point) and
$a'_1\not\in\{b_1,d_1\}$ (since $\bar G$ is loop-free), hence $a'_1=c_1$.
We have $\deg_G(a_1)>3$, $\deg_{\bar G_1}(\bar a_1)<3$, and
$\deg_{\bar G_1}(\bar a_1)=\deg_G(a_1)-2$, hence
$\deg_G(a_1)=4$ and $\deg_{\bar G_1}(\bar a_1)=2$.
This means that the only vertices connected to $a_1$ are
$a$, $c'$, $b_1$, $d_1$. Since $a'_1=c_1$, we have $\deg_G(c_1)=4$,
hence the vertices connected to $c_1$ are
$c$, $a'$, $b_1$, $d_1$. Thus, 
$\sigma(\{a,c',b_1,d_1\})=\{c,a',b_1,d_1\}$. Since $a\mapsto a'$, $c'\mapsto c$,
and $b_1\not\mapsto b_1$, we have $b'_1=d_1$
and we conclude that $G$ is as in Figure \figGfive.
\qed\enddemo

Theorem \thMain\ easily follows from Theorem \thMin\ combined with Theorem 1.
Indeed, By Theorem \thRP, any basic net $\Gamma$ on $\RP^2$
with more than one crossing is dual to a quadrangulation from $\bar\Cal Q$.
Hence $\Gamma$ can be obtained
by successive special face splittings and vertex surroundings
starting from a net dual to a $\bar\Cal Q$-minimal quadrangulation of $\RP^2$.
By Theorem \thMin, $\bar\Cal Q_{\min}$ consists of $\check g_2^5$
and the quotients of those $\Cal Q$-minimal quadrangulations of $S^2$ which admit
a fix point free involution. By Theorem \thNakamoto, $\Cal Q_{\min}=\{W_n\mid n\ge3\}$
(double wheels). It is easy to check that $W_n$ admits a fix point free involution
if and only if $n$ is odd. Thus,
$\bar\Cal Q_{\min}=\{\check g_2^5\}\cup\{\widetilde W_n\mid n$ is odd, $n\ge3\}$.
It remains to note that $g_2^5$ is obtained from $g^1$ by a vertex surrounding
and that the nets $\overline{(2\times n)^*}$ are dual to $\widetilde W_n$.


\head\sectComput. Computations
\endhead

Of course, the best way to generate basic nets in $\RP^2$ is to write
a program based on Theorem \thMain\ and similar to {\tt plantri}
[\refPlantriPap, \refPlantri] or, maybe, just to modify {\tt plantri}.
However, it takes too much efforts for
somebody (like me) who is not familiar with {\tt plantri} internal structure,
so, I used a more lazy approach: I wrote a simple filter {\tt ppf}
for {\tt plantri} (see [\refPPF]).
It reads the output of {\tt plantri} and selects
only those planar graphs which admit an orientation reversing involution
without fix points and fix edges.
Since {\tt plantri} called with {\tt -c2q} option
generates all simple quadrangulations of $S^2$, Theorem \thRP\ ensures
that we obtain in this way all simply embedded quadrangulations of $\RP^2$
(the dual graphs of basic nets with $\ge3$ crossings) without repetitions and omissions.

%
This method is very slow, for example, we need to treat $5.45\cdot 10^{13}$
simple quadrangulations of $S^2$ with $38$ vertices to select only $1735808$
simply embedded quadrangulations of $\RP^2$.
Fortunately, {\tt plantri} is so efficient that this can be done.

The program {\tt ppf} can be used in pipe with {\tt plantri}, for example:
$$
\text{\tt plantri -c2q 18 | ppf}                   \eqno(\eqCommand)
$$
The output is almost the same as the {\tt plantri}'s ascii output but:
\roster
\item"$\bullet$"
       The names of vertices are changed from
       $a,b,c,\dots$ to $a,b,c,\dots, A,B,C,\dots$ so
       that the involution maps $a\mapsto A$, $b\mapsto B$, etc.
\item"$\bullet$"
       We list the neighbourhoods of the lowercase
       vertices only.
\endroster
For example, the first output line produced by the command (\eqCommand) is
$$
\text{\tt 9 bcdef,aDg,agF,aFBH,aHI,aICD,bhic,DEg,gEF}
$$
which corresponds to the net $g^8_1$ in Figure \figBNRP.
The corresponding net in $S^2$ is depicted in Figure \figPPF\
where $S^2$ is supposed to be glued out of the two disks so that
the region names match each other.

\midinsert
\epsfxsize=60mm
\centerline{\epsfbox{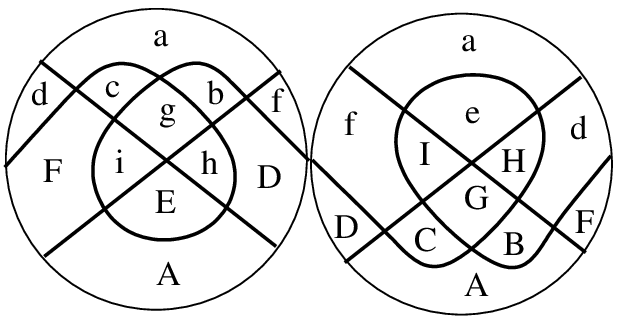}}
\botcaption{ Figure \figPPF }
\endcaption
\endinsert

In Table 1,
$\bar q(\check n)$ is the number of simply embedded
quadrangulations of $\RP^2$ with $\check n=n+1$ vertices (the same
as the number of basic nets on $\RP^2$ with $n$ vertices),
$\bar q_{\bip}(\check n)$ is the number of those of them which are bipartite
(the number of homologically trivial basic nets on $\RP^2$ with $n$ vertices),
and $q_2(2\check n)$ is the number of simple quadrangulations of $S^2$
with $2\check n$ vertices (the same as in [\refBGGMTW; Table 2]), so,
$q_2(2\check n)$ is the number of quadrangulations needed to be checked
in our computation of $\bar q(\check n)$.

\midinsert
\centerline{\hbox{%
\vbox{\offinterlineskip
\halign{&\vrule#&\strut\hfill\;\;\;#\;\;\hfill\cr
\omit& $n$ &         &$\check n$
                           &&$\bar q(\check n)$
                                       &&$\bar q_{\bip}(\check n)$
                                                 &&  $ q_2(2\check n)$
                                                                  &\omit\cr
\omit&\omit&height2pt&\omit&& \omit    && \omit  && \omit         &\omit\cr
\noalign{\hrule}
\omit&\omit&height2pt&\omit&& \omit    && \omit  && \omit         &\omit\cr
\omit&  3  &         &  4  &&    1     &&   0    &&   1           &\omit\cr
\omit&\omit&height2pt&\omit&& \omit    && \omit  && \omit         &\omit\cr
\omit&  4  &         &  5  &&    0     &&   0    &&   1           &\omit\cr
\omit&\omit&height2pt&\omit&& \omit    && \omit  && \omit         &\omit\cr
\omit&  5  &         &  6  &&    2     &&   1    &&   3           &\omit\cr
\omit&\omit&height2pt&\omit&& \omit    && \omit  && \omit         &\omit\cr
\omit&  6  &         &  7  &&    3     &&   2    &&   12          &\omit\cr
\omit&\omit&height2pt&\omit&& \omit    && \omit  && \omit         &\omit\cr
\omit&  7  &         &  8  &&    6     &&   3    &&   64          &\omit\cr
\omit&\omit&height2pt&\omit&& \omit    && \omit  && \omit         &\omit\cr
\omit&  8  &         &  9  &&    12    &&   7    &&   510         &\omit\cr
\omit&\omit&height2pt&\omit&& \omit    && \omit  && \omit         &\omit\cr
\omit&  9  &         &  10 &&    37    &&  22    &&   5146        &\omit\cr
\omit&\omit&height2pt&\omit&& \omit    && \omit  && \omit         &\omit\cr
\omit&  10 &         &  11 &&    95    &&  57    &&   58782       &\omit\cr
\omit&\omit&height2pt&\omit&& \omit    && \omit  && \omit         &\omit\cr
\omit&  11 &         &  12 &&   293    && 174    &&   716607      &\omit\cr
\omit&\omit&height2pt&\omit&& \omit    && \omit  && \omit         &\omit\cr
\omit&  12 &         &  13 &&   923    && 554    &&   9062402     &\omit\cr
\omit&\omit&height2pt&\omit&& \omit    && \omit  && \omit         &\omit\cr
\omit&  13 &         &  14 &&  3086    && 1848   &&   117498072   &\omit\cr
\omit&\omit&height2pt&\omit&& \omit    && \omit  && \omit         &\omit\cr
\omit&  14 &         &  15 && 10504    && 6291   &&   1553048548  &\omit\cr
\omit&\omit&height2pt&\omit&& \omit    && \omit  && \omit         &\omit\cr
\omit&  15 &         &  16 && 36954    && 22052  &&  20858998805  &\omit\cr
\omit&\omit&height2pt&\omit&& \omit    && \omit  && \omit         &\omit\cr
\omit&  16 &         &  17 && 131590   && 78361  && 284057538480  &\omit\cr
\omit&\omit&height2pt&\omit&& \omit    && \omit  && \omit         &\omit\cr
\omit&  17 &         &  18 && 475793   && 282420 && 3915683667721 &\omit\cr
\omit&\omit&height2pt&\omit&& \omit    && \omit  && \omit         &\omit\cr
\omit&  18 &         &  19 &&  1735808  && 1027336  &&  54565824458485 &\omit\cr
\omit&\omit&height2pt&\omit&& \omit    && \omit  && \omit         &\omit\cr
}\hrule}}}
\botcaption{Table 1}
      Basic nets on $\RP^2$ with $n$ and on $S^2$ with $2n$ vertices 
\endcaption
\endinsert

\Refs
\def\r{\ref}
\r\no\refBGGMTW
\by    G.~Brinkmann, S.~Greenberg, C.~Greenhill, B.~D.~McKay, R.~Thomas, P.~Wollan
\paper Generation of simple quadrangulations of the sphere
\jour  Discrete Math. \vol 305 \yr 2005 \pages 33--54
\endref

\r\no\refPlantriPap
\by     G.~Brinkmann, B.~D.~McKay
\paper  Fast generation of planar graphs
\jour   MATCH: Commun. Math. Comput. Chem. \vol 58\yr 2007 \pages 323--357
\transl Expanded:
        http://cs.anu.edu.au/$\tilde{\;\,}$bdm/papers/plantri-full.pdf
\endref

\r\no\refPlantri
\by    G.~Brinkmann, B.~D.~McKay
\paper The program plantri 
\jour      http://cs.anu.edu.au/$\tilde{\;\,}$bdm/plantri
\endref

\r\no\refC
\by     J.~H.~Conway
\paper  An enumeration of knots and links, and some of their algebraic
        properties
\inbook Computational Problems of Abstract Algebra (Proc. Conf., Oxford, 1967)
\publ   Pergamon \publaddr Oxford \yr 1970 \pages 329--358
\transl Available at http://www.math.ed.ac.uk/$\tilde{\;}$aar/knots/conway.pdf
\endref


\r\no\refD
\by    J.~Drobotukhina
\paper Classification of links in $\Bbb RP^3$ with at most six crossings
\inbook in ``Topology of manifolds and varieties'' (ed. O.Ya.~Viro)
\bookinfo Advances in Soviet Math. \vol 18 \yr 1994 \publ A.M.S. 
\publaddr \pages 87--121
\endref

\r\no\refN
\by    A.~Nakamoto
\paper Generating Quadrangulations of Surfaces with Minimum Degree at Least 3
\jour  Journal of Graph Theory \vol 30 \yr 1999 \issue 3 \pages 223-234
\endref

\r\no\refPPF
\by    S.~Yu.~Orevkov
\paper The program ppf 
\jour  http://picard.ups-tlse.fr/$\tilde{\;}$orevkov/ppf.c
\endref

\endRefs
\enddocument